# Characterization of non-differentiable points in a function by Fractional derivative of Jumarrie type

Uttam Ghosh (1), Srijan Sengupta(2), Susmita Sarkar (2), Shantanu Das (3)

(1): Department of Mathematics, Nabadwip Vidyasagar College, Nabadwip, Nadia, West Bengal, India;
Email: uttam_math@yahoo.co.in
(2):Department of Applied Mathematics, Calcutta University, Kolkata, India
Email: susmita62@yahoo.co.in
(3)Scientist H+, RCSDS, BARC Mumbai India
Senior Research Professor, Dept. of Physics, Jadavpur University Kolkata
Adjunct Professor. DIAT-Pune
Ex-UGC Visiting Fellow Dept. of Applied Mathematics, Calcutta University, Kolkata India
Email (3): shantanu@barc.gov.in

The Birth of fractional calculus from the question raised in the year 1695 by Marquis de L'Hopital to Gottfried Wilhelm Leibniz, which sought the meaning of Leibniz's notation for the derivative of order N when N = 1/2. Leibnitz responses it is an apparent paradox from which one day useful consequences will be drown.

## Abstract

There are many functions which are continuous everywhere but not differentiable at some points, like in physical systems of ECG, EEG plots, and cracks pattern and for several other phenomena. Using classical calculus those functions cannot be characterized-especially at the non-differentiable points. To characterize those functions the concept of Fractional Derivative is used. From the analysis it is established that though those functions are unreachable at the non-differentiable points, in classical sense but can be characterized using Fractional derivative. In this paper we demonstrate use of modified Riemann-Liouvelli derivative by Jumarrie to calculate the fractional derivatives of the non-differentiable points of a function, which may be one step to characterize and distinguish and compare several non-differentiable points in a system or across the systems. This method we are extending to differentiate various ECG graphs by quantification of non-differentiable points; is useful method in differential diagnostic. Each steps of calculating these fractional derivatives is elaborated.

## 1.0   Introduction

The concept of classical calculus in modern form was developed in end seventeenth by Newton and Leibnitz [1]. Leibniz used the symbol $\frac{d^n y}{dx^n}$ to denote the *n*-th order derivative of $f(x)$ [2]. From the above developed notation de L'Hospital asked Leibniz what is the meaning of $\frac{d^n y}{dx^n}$, for $n = 1/2$ gives the birth of fractional derivative. Recently authors [3-8] are trying to generalize the concept of derivative for all real and complex values of *n*. Again in generalized notation when *n* is positive it will be derivative and for negative *n* it will be the notion of integration. In these early methods the derivative of constant is found non-zero. So $\frac{d^n y}{dx^n}$ $\forall n$ denotes the generalized order derivative and integration, or generalized differ-integration. The basic definition of generalized derivative are the formulas from Grunwald-Letinikov(G-L) definition, Riemann-Liouville (R-L) and Caputo definition [6,8]. The Riemann-Liouvelli definition retuns a non-zero for fractional derivative of a constant. This differs from the basic definition of classical derivative. To overcome this gap Jumarie [11] modified the definition of the fractional order derivative of left Riemann-Liouville. In this paper we have modified the right Riemann-Liouville fractional derivative and used both the modified definition (left and right) of derivative to find the derivative of the non-differentiable functions and the result is interpreted graphically. The calculus on rough unreachable functions is developed via Kolwankar-Gangal (KG) derivative a Local Fractional Derivative (LFD) [16]-[20], where limit is taken at unreachable point to get LFD. In this paper here there is no limit concept; instead the modification is done on classical Rieman-Liouvelli fractional derivative by constructing an offset function and doing integration for the defined interval; for the left and right fractional derivative. The organization of the paper is as follows in section 1.1 some definition with examples is given, in section 2.0 fractional derivative of some non-differentiable is calculated with graphical results presented.

**1.1 Some definitions**

**1.11 Grunwald-Letinikov definition**

Let $f(t)$ be any function then the $\alpha$-th order derivative $\alpha \in \mathbb{R}$ of $f(t)$ is defined by

$$\begin{aligned} {}_a D_t^\alpha f(t) &= \lim_{\substack{h \to 0 \\ nh \to t-a}} h^{-\alpha} \sum_{r=0}^{n} \binom{\alpha}{r} f(t - rh) \\ &= \lim_{\substack{h \to 0 \\ nh \to t-a}} h^{-\alpha} \sum_{r=0}^{n} \frac{\alpha!}{r!(\alpha - r)!} f(t - rh) \\ &= \frac{1}{(-\alpha - 1)!} \int_a^t (t - \tau)^{-\alpha - 1} f(\tau) d\tau \\ &= \frac{1}{\Gamma(-\alpha)} \int_a^t \frac{f(\tau)}{(t - \tau)^{\alpha + 1}} d\tau \end{aligned} \quad (1)$$

Where $\alpha$ is any arbitrary number real or complex and $\binom{\alpha}{r} = \frac{\alpha!}{r!(\alpha - r)!} = \frac{\Gamma(\alpha + 1)}{\Gamma(r+1)\Gamma(\alpha - r + 1)}$

The above formula becomes fractional order integration if we replace $\alpha$ by $-\alpha$ which is

$$_aD_t^{-\alpha}f(t) = \frac{1}{\Gamma(\alpha)}\int_a^t (t-\tau)^{\alpha-1}f(\tau)d\tau \qquad (2)$$

Using the above formula we get for $f(t) = (t-a)^\gamma$,

$$_aD_t^{\alpha}(t-a)^\gamma = \frac{1}{\Gamma(-\alpha)}\int_a^t (t-\tau)^{-(\alpha+1)}(\tau-a)^\gamma d\tau$$

Using the substitution $\tau = a + \xi(t-a)$ we have for $\tau = a$, $\xi = 0$ and for $\tau = t, \xi = 1$; $d\tau = (t-a)d\xi$, $(t-\tau) = t-a-\xi(t-a) = (t-a)(1-\xi)$; $(\tau-a) = \xi(t-a)$, we get the following

$$_aD_t^{\alpha}(t-a)^\gamma = \frac{1}{\Gamma(-\alpha)}\int_a^t (t-\tau)^{-(\alpha+1)}(\tau-a)^\gamma d\tau$$

$$= \frac{1}{\Gamma(-\alpha)}\int_0^1 (t-a)^{-\alpha-1}(1-\xi)^{-\alpha-1}\xi^\gamma (t-a)^\gamma (t-a)d\xi$$

$$= \frac{(t-a)^{\gamma-\alpha}}{\Gamma(-\alpha)}\int_0^1 \xi^\gamma (1-\xi)^{-(\alpha+1)}d\xi$$

$$= \frac{(t-a)^{\gamma-\alpha}}{\Gamma(-\alpha)}\mathrm{B}(-\alpha,\gamma+1)$$

$$= \frac{\Gamma(\gamma+1)}{\Gamma(\gamma+1-\alpha)}(t-a)^{\gamma-\alpha}, (\alpha < 0, \gamma > -1)$$

We used Beta-function $\mathrm{B}(-\alpha, \gamma+1) = \int_0^1 \xi^\gamma (1-\xi)^{-(\alpha+1)} d\xi = \frac{\Gamma(-\alpha)\Gamma(\gamma+1)}{\Gamma(-\alpha+\gamma+1)}$ defined as

$$\mathrm{B}(p,q) \overset{\text{def}}{=} \int_0^1 u^{p-1}(1-u)^{q-1} du = \frac{\Gamma(p)\Gamma(q)}{\Gamma(p+q)}$$

If the function $f(t)$ be such that $f^k(t)$, $k = 1,2,3,...,m+1$ is continuous in the closed interval $[a,t]$ and $m \leq \alpha < m+1$ then

$$_aD_t^\alpha f(t) = \sum_{k=0}^m \frac{f^{(k)}(a)(t-a)^{-\alpha+k}}{\Gamma(-\alpha+k+1)} + \frac{1}{\Gamma(-\alpha+m+1)}\int_a^t (t-\tau)^{m-\alpha}f^{(m+1)}(\tau)d\tau \qquad (3)$$

### 1.12 Riemann-Liouville (R-L) definition of fractional derivative

Let the function $f(t)$ is one time integrable then the integro-differential expression

$$_aD_t^\alpha f(t) = \frac{1}{\Gamma(-\alpha+m+1)}\left(\frac{d}{dt}\right)^{m+1}\int_a^t (t-\tau)^{m-\alpha} f(\tau)d\tau \tag{4}$$

is known as the Riemann-Liouville definition of fractional derivative [6] with $m \leq \alpha < m+1$.

In Riemann-Liouville definition the function $f(t)$ is getting fractionally integrated and then differentiated $m+1$ whole-times but in obtaining the formula (3) $f(t)$ must be $m+1$ time differentiable. If the function $f(t)$ is $m+1$ whole times differentiable then the definition (1), (3) and (4) are equivalent.

Using integration by parts formula in (2), that is $_aD_t^{-\alpha}f(t) = \frac{1}{\Gamma(\alpha)}\int_a^t (t-\tau)^{\alpha-1} f(\tau)d\tau$ we get

$$_aD_t^{-\alpha}f(t) = \frac{f(a)(t-a)^\alpha}{\Gamma(\alpha+1)} + \frac{1}{\Gamma(\alpha+1)}\int_a^t (t-\tau)^\alpha f'(\tau)d\tau \tag{5}$$

The left R-L fractional derivative is defined by

$$_aD_t^\alpha f(t) = \frac{1}{\Gamma(k-\alpha)}\left(\frac{d}{dt}\right)^k \int_a^t (t-\tau)^{k-\alpha-1} f(\tau)d\tau \qquad k-1 \leq \alpha < k \tag{6}$$

And the right R-L derivative is

$$_tD_b^\alpha f(t) = \frac{1}{\Gamma(k-\alpha)}\left(-\frac{d}{dt}\right)^k \int_t^b (\tau-t)^{k-\alpha-1} f(\tau)d\tau \qquad k-1 \leq \alpha < k \tag{7}$$

In above definitions $k \in \mathbb{Z}$ that is integer just greater than alpha and $\alpha > 0$, $\alpha \in \mathbb{R}$

From the above definition it is clear that if at time $t$ the function $f(t)$ describes a certain dynamical system developing with time then for $\tau < t$, where $t$ is the present time then state $f(t)$ represent the past time and similarly if $\tau > t$ then $f(t)$ represent the future time. Therefore the left derivative represents the past state of the process and the right hand derivative represents the future stage.

### 1.13 Caputo definition of fractional derivative

In the R-L type definition the initial conditions contains the limit of R-L fractional derivative such as $\lim_{x \to a} {_aD_t^{\alpha-1}} = b_1$ etc that is fractional initial staes. But if the initial conditions are $f(a) = b_1, f'(a) = b_2 \cdots$ type then R-L definition fails and to overcome these problems M. Caputo [15] proposed new definition of fractional derivative in the following form

$$_a^C D_t^\alpha f(t) = \frac{1}{\Gamma(\alpha-n)} \int_a^t \frac{f^{(n)}(\tau)}{(t-\tau)^{\alpha+1-n}} d\tau, \qquad n-1 < \alpha < n \tag{8}$$

Under natural condition on the function $f(t)$ and as $\alpha \to n$ the Caputo derivative becomes a conventional *n*-th order derivative of the function. The main advantage of the Caputo derivative is the initial conditions of the fractional order derivatives are conventional derivative type-i.e requiring integer order states.

In R-L derivative the derivative of Constant (C) is non-zero. Since

$$_0 D_t^{-\alpha} C = \frac{1}{\Gamma(-\alpha)} \int_0^t (t-\tau)^{-\alpha-1} C d\tau$$

$$= \frac{C t^{-\alpha}}{\Gamma(-\alpha)} \int_0^1 (1-x)^{-\alpha-1} dx, \qquad \text{Where } \tau = tx$$

$$= \frac{C t^{-\alpha}}{\Gamma(-\alpha)} \left[ \frac{(1-x)^{-\alpha}}{-\alpha} \right]_0^1$$

$$= \begin{cases} \dfrac{C t^{-\alpha}}{\Gamma(1-\alpha)} & \text{for } \alpha < 0 \\ 0 & \text{for } \alpha > 0 \end{cases}$$

## 1.14 Jumarie definition of fractional derivative

On the other hand to overcome the misconception derivative of a constant is zero in the conventional integer order derivative Jumarie [11] revised the R-L derivative in the following form

$$D_x^\alpha f(x) = \frac{1}{\Gamma(-\alpha)} \int_0^x (x-\xi)^{-\alpha-1} f(\xi) d\xi, \qquad \text{for} \qquad \alpha < 0$$

$$= \frac{1}{\Gamma(1-\alpha)} \frac{d}{dx} \int_0^x (x-\xi)^{-\alpha} [f(\xi) - f(0)] d\xi, \qquad \text{for} \qquad 0 < \alpha < 1$$

$$= \left( f^{(\alpha-n)}(x) \right)^{(n)} \qquad \text{for} \qquad n \leq \alpha < n+1, \qquad n \geq 1.$$

The above definition [11] is developed using left R-L derivative. Similarly using the right R-L derivative other type can be develop. Note in the above definition for negative fractional orders the expression is just Riemann-Liouvelli fractional integration. The modification is carried out in R-L the derivative formula, for the positive fractional orders alpha. The idea is to remove the offset value of function at start point of the fractional derivative from the function, and carry out R-L derivative usually done for the function.

First we want to find the derivative of constant (C) using right R-L derivative,

$$D_t^\alpha(C) = \frac{1}{\Gamma(1-\alpha)}\left(-\frac{d}{dx}\right)\int_x^b (\xi-x)^{-\alpha} C\, d\xi$$

$$= \begin{cases} \dfrac{C}{\Gamma(1-\alpha)}(b-x)^{-\alpha}, & \text{for} \quad \alpha < 0 \\ 0, & \text{other wise} \quad \text{for} \quad (\alpha > 0) \end{cases}$$

Since for any function $f(x)$ in the interval $[a,b]$ which satisfies the conditions of modified fractional derivative [6] can be written as

$$f(x) = f(b) - [f(b) - f(x)]$$
$$D^\alpha f(x) = D^\alpha f(b) - D^\alpha [f(b) - f(x)]$$
For $\quad \alpha < 0,$

$$D_x^\alpha [f(b) - f(x)] = -\frac{1}{\Gamma(-\alpha)}\int_x^b (\xi-x)^{-\alpha-1} f(\xi)\, d\xi$$

For $\quad 0 < \alpha < 1,$

$$D_x^\alpha [f(b) - f(x)] = -f^{(\alpha)}(x) = -\left(f^{(\alpha-1)}(x)\right)'$$

$$f^{(\alpha)}(x) = -\frac{1}{\Gamma(1-\alpha)}\frac{d}{dx}\int_x^b (\xi-x)^{-\alpha}[f(b)-f(\xi)]\, d\xi$$

and for $n \leq \alpha < n+1,$
$$f^{(\alpha)}(x) = [f^{(\alpha-n)}(x)]^{(n)}$$

Thus finally we can define in the following form

$$D_x^\alpha f(x) = -\frac{1}{\Gamma(-\alpha)}\int_x^b (\xi-x)^{-\alpha-1} f(\xi)\, d\xi, \quad \text{for} \quad \alpha < 0$$

$$= -\frac{1}{\Gamma(1-\alpha)}\frac{d}{dx}\int_x^b (\xi-x)^{-\alpha}[f(b)-f(\xi)]\, d\xi, \quad \text{for} \quad 0 < \alpha < 1$$

$$= \left(f^{(\alpha-n)}(x)\right)^{(n)} \quad \text{for} \quad n \leq \alpha < n+1, \quad n \geq 1.$$

**Example**: Use the above definition to a continuous and differentiable function

$$f(x) = x - c, \quad a \leq x \leq b.$$

By using the Jumarie definition (as above) we obtain for left fractional derivative

$$f_L^{(\alpha)}(x) = \frac{1}{\Gamma(1-\alpha)} \frac{d}{dx} \int_a^x (x-\xi)^{-\alpha} [f(\xi) - f(a)] d\xi, \qquad 0 < \alpha < 1.$$

$$= \frac{1}{\Gamma(1-\alpha)} \frac{d}{dx} \int_a^x (x-\xi)^{-\alpha} [\xi - c - (a-c)] d\xi = \frac{1}{\Gamma(1-\alpha)} \frac{d}{dx} \int_a^x (x-\xi)^{-\alpha} [\xi - a] d\xi$$

$$= \frac{1}{\Gamma(1-\alpha)} \frac{d}{dx} \int_a^x (x-\xi)^{-\alpha} [-(x-\xi) + (x-a)] d\xi$$

$$= \frac{1}{\Gamma(1-\alpha)} \frac{d}{dx} \int_a^x [-(x-\xi)^{1-\alpha} + (x-a)(x-\xi)^{-\alpha}] d\xi$$

$$= \frac{1}{\Gamma(1-\alpha)} \frac{d}{dx} \left\{ -\left[ -\frac{(x-\xi)^{2-\alpha}}{2-\alpha} \right]_a^x + (x-a)\left[ -\frac{(x-\xi)^{1-\alpha}}{1-\alpha} \right]_a^x \right\}$$

$$= \frac{-1}{\Gamma(1-\alpha)} \frac{d}{dx} \left[ \frac{(x-a)^{2-\alpha}}{2-\alpha} - \frac{(x-a)^{2-\alpha}}{1-\alpha} \right]$$

$$= \frac{1}{\Gamma(1-\alpha)} \frac{d}{dx} \left[ \frac{(x-a)^{2-\alpha}}{(2-\alpha)(1-\alpha)} \right] = \frac{(x-a)^{1-\alpha}}{\Gamma(2-\alpha)} \qquad \text{we used here} \qquad n\Gamma(n) = \Gamma(n+1)$$

Therefore

$$f_L^{(\alpha)}\left(\frac{a+b}{2}\right) = \frac{1}{\Gamma(2-\alpha)} \left(\frac{a+b}{2} - a\right)^{1-\alpha} = \frac{\left(\frac{b-a}{2}\right)^{1-\alpha}}{\Gamma(2-\alpha)}$$

Again using our right modified definition we obtain

$$f_R^{(\alpha)}(x) = \frac{1}{\Gamma(1-\alpha)} \left(-\frac{d}{dx}\right) \int_x^b (\xi - x)^{-\alpha} [f(b) - f(\xi)] d\xi, \qquad 0 < \alpha < 1.$$

$$= \frac{1}{\Gamma(1-\alpha)} \frac{d}{dx} \int_x^b (\xi - x)^{-\alpha} [\xi - c - (b-c)] d\xi$$

$$= \frac{1}{\Gamma(1-\alpha)} \frac{d}{dx} \int_x^b (\xi - x)^{-\alpha} [(\xi - x) + (x - b)] d\xi$$

$$= \frac{1}{\Gamma(1-\alpha)} \frac{d}{dx} \int_x^b [(\xi - x)^{1-\alpha} + (x - b)(\xi - x)^{-\alpha}] d\xi$$

$$= \frac{1}{\Gamma(1-\alpha)} \frac{d}{dx} \left[ \frac{(\xi-x)^{2-\alpha}}{2-\alpha} + (x-b) \frac{(\xi-x)^{1-\alpha}}{1-\alpha} \right]_x^b$$

$$= \frac{1}{\Gamma(1-\alpha)} \frac{d}{dx} \left[ \frac{(b-x)^{2-\alpha}}{2-\alpha} - \frac{(b-x)^{2-\alpha}}{1-\alpha} \right]$$

$$= \frac{(b-x)^{1-\alpha}}{\Gamma(2-\alpha)}$$

Therefore

$$f_R^{(\alpha)}\left(\frac{a+b}{2}\right) = \frac{1}{\Gamma(2-\alpha)} \left(b - \frac{a+b}{2}\right)^{1-\alpha} = \frac{\left(\frac{b-a}{2}\right)^{1-\alpha}}{\Gamma(2-\alpha)}$$

Thus in both the cases (for *L* and *R*) value of $f^{(\alpha)}\left(\frac{a+b}{2}\right)$ is equal. Thus for continuous and differentiable functions both the values are equal, and is equal to $\dfrac{\left(\frac{b-a}{2}\right)^{1-\alpha}}{\Gamma(2-\alpha)}$

## 1.15 Unreachable function

There are many functions which are continuous for all *x* but not-differentiable at some points or at all points. These functions are named as **unreachable functions**. The function $f(x) = |x - 1/2|$ is unreachable at the point $x = 1/2$. The function $f(x) = \tan(x)$ is unreachable at the point $x = \pi/2$.

These functions are non-differentiable i.e. unreachable functions, at some points in the interval. To study those functions we are considering the following examples and the fraction derivative will be found using Jumarie modified definition.

## 2.0 Fractional derivative of Some Unreachable Functions

**Example 1**: $f(x) = |x - (1/2)|$ defined on [0, 1]. This function is continuous for all *x* in the given interval but not differentiable at $x = 1/2$ i.e. we cannot draw tangent at this point to the curve. The curve is symmetric about the non-differentiable point which is clear from the figure 1. To study behavior of the function at we want to find out fractional derivative at $x = 1/2$ using the modified fractional derivative of Jumarie.

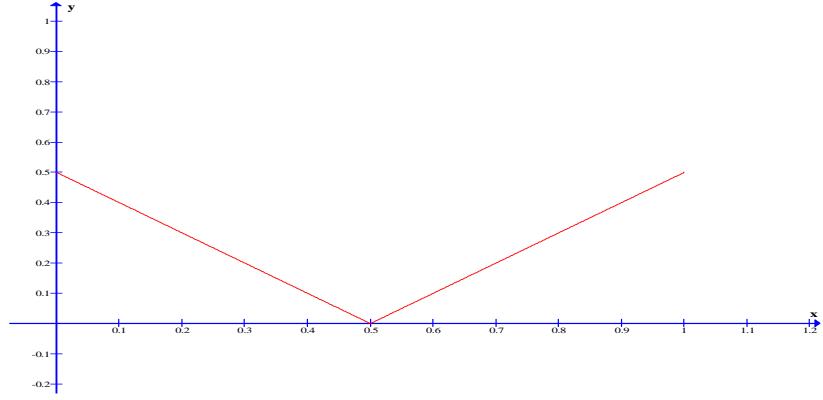

**Fig. 1.0 Graph of the function |x-1/2|**

(a) The fractional order derivative using Jumarie modified definition is

$$f_L^{(\alpha)}(x) = \frac{1}{\Gamma(1-\alpha)} \frac{d}{dx} \int_0^x (x-\xi)^{-\alpha} [f(\xi) - f(0)] d\xi, \qquad 0 < \alpha < 1.$$

When $0 \leq x \leq 1/2$ $\quad f(x) = -x + (1/2)$

$$f_L^{(\alpha)}(x) = \frac{1}{\Gamma(1-\alpha)} \frac{d}{dx} \int_0^x (x-\xi)^{-\alpha} [f(\xi) - f(0)] d\xi = \frac{1}{\Gamma(1-\alpha)} \frac{d}{dx} \int_0^x (x-\xi)^{-\alpha} [-\xi + (1/2) - (1/2)] d\xi$$

$$= \frac{1}{\Gamma(1-\alpha)} \frac{d}{dx} \int_0^x (x-\xi)^{-\alpha} [(x-\xi) - x] d\xi = \frac{1}{\Gamma(1-\alpha)} \frac{d}{dx} \int_0^x [(x-\xi)^{1-\alpha} - x(x-\xi)^{-\alpha}] d\xi$$

$$= \frac{1}{\Gamma(1-\alpha)} \frac{d}{dx} \left[ -\frac{(x-\xi)^{2-\alpha}}{2-\alpha} + x \frac{(x-\xi)^{1-\alpha}}{1-\alpha} \right]_0^x$$

$$= \frac{1}{\Gamma(1-\alpha)} \frac{d}{dx} \left[ -\frac{x^{2-\alpha}}{(2-\alpha)(1-\alpha)} \right] = -\frac{x^{1-\alpha}}{\Gamma(2-\alpha)}$$

Again when $1/2 \leq x \leq 1$ $\quad f(x) = x - (1/2)$ the fractional derivative from zero to half and beyond is done in two separate segments as below

$$f_L^{(\alpha)}(x) = \frac{1}{\Gamma(1-\alpha)} \frac{d}{dx} \left( \int_0^{1/2} + \int_{1/2}^x \right)(x-\xi)^{-\alpha}[f(\xi)-f(0)]d\xi$$

$$= \frac{1}{\Gamma(1-\alpha)} \frac{d}{dx}\left[ \int_0^{1/2}(x-\xi)^{-\alpha}[-\xi+(1/2)-(1/2)]d\xi + \int_{1/2}^x (x-\xi)^{-\alpha}[\xi-(1/2)-(1/2)]d\xi \right]$$

$$= \frac{1}{\Gamma(1-\alpha)} \frac{d}{dx}\left[ \left(-\frac{(x-\xi)^{2-\alpha}}{2-\alpha} + x\frac{(x-\xi)^{1-\alpha}}{1-\alpha}\right)_0^{1/2} + \int_{1/2}^x (x-\xi)^{-\alpha}[\xi-1]d\xi \right]$$

$$= \frac{1}{\Gamma(1-\alpha)} \frac{d}{dx}\left[ \left(-\frac{(x-\xi)^{2-\alpha}}{2-\alpha} + x\frac{(x-\xi)^{1-\alpha}}{1-\alpha}\right)_0^{1/2} + \int_{1/2}^x (x-\xi)^{-\alpha}[(x-1)-(x-\xi)]d\xi \right]$$

$$= \frac{1}{\Gamma(1-\alpha)} \frac{d}{dx}\left[ \left(-\frac{(x-\xi)^{2-\alpha}}{2-\alpha} + x\frac{(x-\xi)^{1-\alpha}}{1-\alpha}\right)_0^{1/2} - \left(-\frac{(x-\xi)^{2-\alpha}}{2-\alpha} + (x-1)\frac{(x-\xi)^{1-\alpha}}{1-\alpha}\right)_{1/2}^x - \right]$$

$$= \frac{1}{\Gamma(1-\alpha)} \frac{d}{dx}\left[ -\frac{(x-1/2)^{2-\alpha}-x^{2-\alpha}}{2-\alpha} + x\frac{(x-1/2)^{1-\alpha}-x^{1-\alpha}}{1-\alpha} - \frac{(x-1/2)^{2-\alpha}}{2-\alpha} + (x-1)\frac{(x-1/2)^{1-\alpha}}{1-\alpha} \right]$$

$$= \frac{2(x-1/2)^{1-\alpha} - x^{1-\alpha}}{\Gamma(2-\alpha)}$$

Therefore

$$f_L^{(\alpha)}(x) = \begin{cases} -\dfrac{x^{1-\alpha}}{\Gamma(2-\alpha)}, & 0 \leq x \leq 1/2 \\ \dfrac{2(x-1/2)^{1-\alpha} - x^{1-\alpha}}{\Gamma(2-\alpha)}, & 1/2 \leq x \leq 1 \end{cases}$$

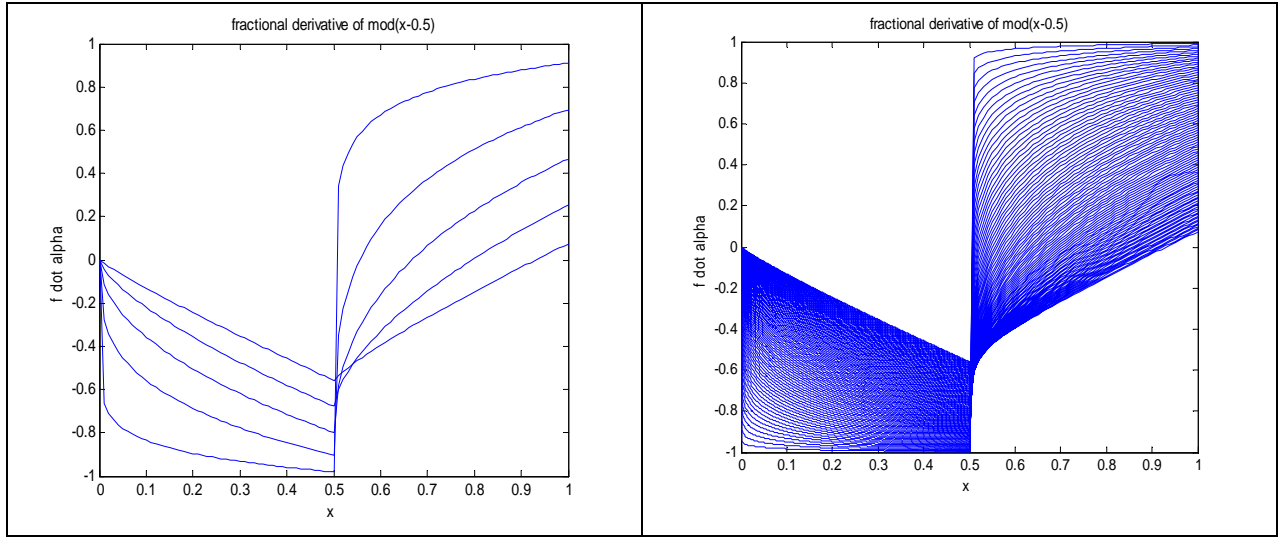

**Fig.2 Graph of the function** $f_L^{(\alpha)}(x)$ **for different values of alpha.**

From the above expression it is clear that though $f'(1/2)$ does not exists but the fractional derivative $f_L^{(\alpha)}$ exists at $x = 1/2$ and equals to $f_L^{(\alpha)}(1/2) = -\dfrac{(1/2)^{2-\alpha}}{\Gamma(2-\alpha)}$. For $\alpha = 0.5$ we get the value as $f_L^{(0.5)}(1/2) = -\dfrac{(1/2)^{3/2}}{\Gamma(3/2)}$ with $\Gamma(3/2) = \sqrt{\pi}/2$ we get a value of left half derivative at $x = 1/2$ as $f_L^{(0.5)}(1/2) = -1\big/\sqrt{2\pi}$

(b) The fractional order derivative using right R-L definition and modifying the same as Jumarrie, on same function we get

$$f_R^{(\alpha)}(x) = -\frac{1}{\Gamma(1-\alpha)}\frac{d}{dx}\int_x^1 (\xi - x)^{-\alpha}[f(1) - f(\xi)]d\xi, \qquad 0 < \alpha < 1.$$

When $0 \leq x \leq 1/2 \qquad f(x) = -x + (1/2)$

$$f_R^{(\alpha)}(x) = -\frac{1}{\Gamma(1-\alpha)}\frac{d}{dx}\left(\int_x^{1/2} + \int_{1/2}^1\right)(\xi - x)^{-\alpha}[f(1) - f(\xi)]d\xi$$

$$= \frac{1}{\Gamma(1-\alpha)}\frac{d}{dx}\left(\int_x^{1/2}(\xi - x)^{-\alpha}[-x - (\xi - x)]d\xi + \int_{1/2}^1 [(\xi - x)^{1-\alpha} + (x-1)(\xi - x)^{-\alpha}]d\xi\right)$$

$$= \frac{1}{\Gamma(1-\alpha)}\frac{d}{dx}\left[\begin{array}{l}\dfrac{-x(1/2-x)^{1-\alpha}}{1-\alpha} - \dfrac{(1/2-x)^{2-\alpha}}{2-\alpha} + \dfrac{(1-x)^{2-\alpha} - (1/2-x)^{2-\alpha}}{2-\alpha} \\ +(x-1)\dfrac{(1-x)^{1-\alpha} - (1/2-x)^{1-\alpha}}{1-\alpha}\end{array}\right]$$

$$= \frac{(1-x)^{1-\alpha}}{\Gamma(2-\alpha)}$$

When $1/2 \leq x \leq 1 \qquad f(x) = x - (1/2)$

$$f_R^{(\alpha)}(x) = \frac{1}{\Gamma(1-\alpha)} \frac{d}{dx} \int_x^1 (\xi - x)^{-\alpha} (\xi - (1/2) - (1/2))] d\xi$$

$$= \frac{1}{\Gamma(1-\alpha)} \frac{d}{dx} \int_x^1 (\xi - x)^{-\alpha} [(\xi - x) + (x - 1)] d\xi$$

$$= \frac{1}{\Gamma(1-\alpha)} \frac{d}{dx} \left[ \frac{(\xi - x)^{2-\alpha}}{2-\alpha} + (x-1) \frac{(\xi - x)^{1-\alpha}}{1-\alpha} \right]_x^1$$

$$= \frac{1}{\Gamma(1-\alpha)} \frac{d}{dx} \left[ \frac{(1-x)^{2-\alpha}}{2-\alpha} + (x-1) \frac{(1-x)^{1-\alpha}}{1-\alpha} \right]$$

$$= \frac{(1-x)^{1-\alpha}}{\Gamma(2-\alpha)}$$

Therefore

$$f_R^{(\alpha)}(x) = \begin{cases} \dfrac{\{(1-x)^{1-\alpha} - 2(1/2-x)^{1-\alpha}\}}{\Gamma(2-\alpha)}, & \text{for} \quad 0 \leq x \leq 1/2 \\ \dfrac{(1-x)^{1-\alpha}}{\Gamma(2-\alpha)}, & \text{for} \quad 1/2 \leq x \leq 1 \end{cases}$$

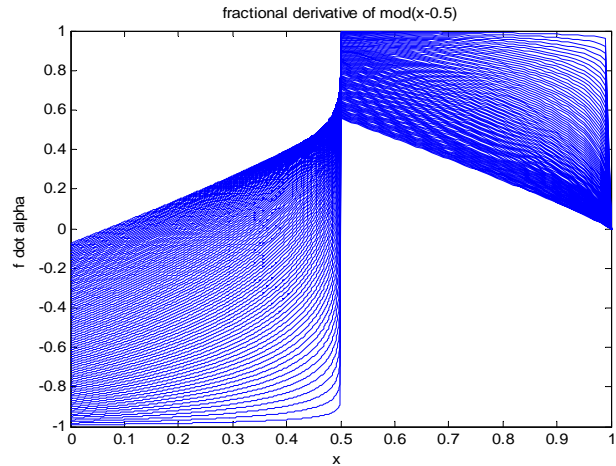

**Fig.3 Graph of the function $f_R^{(\alpha)}(x)$ for different values of alpha.**

From figure 3 it is clear that the right modified derivative exist for this non-differentiable function. Thus both the cases we noticed that the function is not differentiable at $x = 1/2$ but its fractional order derivative exists. The value of half right derivative $f_R^{(0.5)}(1/2) = 1/\sqrt{2\pi}$

From the above two examples it is clear that for differentiable functions the modified definition (both left and right) of the fractional derivative gives the same value at any particular point but for those functions having non-differentiability at some point gives different value for the **LFFT** and RIGHT MODIFIED DERIVATIVE. The difference in values of the fractional derivative at the non-differentiable points indicates the **PHASE transition** at the non-differentiable points. **The difference of the LFFT and RIGHT MODIFIED DERIVATIVE is here defining as the indicator of level of phase transition.**

In example-1 we consider a function which is symmetric about the non-differentiable point and the function is linear in both sides of the non-differentiable points. Now we are considering a function which is non-symmetric with respect to the non-differentiable point and linear in both side of the non-differentiable point.

**Example 2**: Let

$$f(x) = \begin{cases} 10x - 16, & 2 \leq x \leq 2.5 \\ 49 - 16x, & 2.5 \leq x \leq 3 \end{cases}$$

which arises in approximation of a lead in ECG graph of V5 peak of a patient

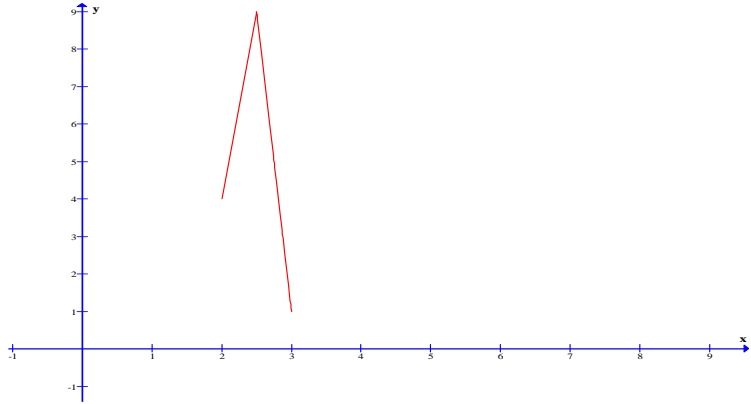

**Fig-4: ECG peak of V5**

From the figure it is clear that this function is continuous for all values of $x$ in the given interval and non-symmetric about the point $x = 2.5$ (the non-differentiable point). We give a translation $z = x - 2$ and rewrite the functional form by same notation and the interval $[2,3]$ to $[0,1]$ and the translated function is

$$f(x) = \begin{cases} 10x + 4, & 0 \leq x \leq 0.5 \\ 17 - 16x, & 0.5 \leq x \leq 1 \end{cases}.$$

This function is continuous at for all $x$ but not-differentiable at $x = 1/2$ but nature of discontinuity is different from the function $f(x) = \left| x - \frac{1}{2} \right|$ in [0, 1].

(a) The fractional order derivative using Jumarie modified definition is

$$f_L^{(\alpha)}(x) = \frac{1}{\Gamma(1-\alpha)} \frac{d}{dx} \int_0^x (x-\xi)^{-\alpha} [f(\xi) - f(0)] d\xi, 0 < \alpha < 1$$

$$= \frac{1}{\Gamma(1-\alpha)} \frac{d}{dx} \int_0^x (x-\xi)^{-\alpha} [f(\xi) - 4] d\xi$$

When $0 \leq x \leq 1/2$     $f(x) = 10x + 4$     $f(0) = 4$

$$f_L^{(\alpha)}(x) = \frac{1}{\Gamma(1-\alpha)} \frac{d}{dx} \int_0^x (x-\xi)^{-\alpha} [10\xi + 4 - 4] d\xi = \frac{10}{\Gamma(1-\alpha)} \frac{d}{dx} \int_0^x [x(x-\xi)^{-\alpha} - (x-\xi)^{1-\alpha}] d\xi$$

$$= \frac{10}{\Gamma(1-\alpha)} \frac{d}{dx} \left[ \frac{(x-\xi)^{2-\alpha}}{2-\alpha} - x\frac{(x-\xi)^{1-\alpha}}{1-\alpha} \right]_0^x = \frac{10}{\Gamma(1-\alpha)} \frac{d}{dx} \left[ \frac{x^{2-\alpha}}{1-\alpha} - \frac{x^{2-\alpha}}{2-\alpha} \right]$$

$$= \frac{10}{\Gamma(1-\alpha)} \frac{d}{dx} \left[ \frac{x^{2-\alpha}}{(2-\alpha)(1-\alpha)} \right] = \frac{10 x^{1-\alpha}}{\Gamma(2-\alpha)}$$

When $0 \geq x \geq 1/2$     $f(x) = 17 - 16x$ here we require the value at $f(0) = 4$ and also interval $[0, 0.5]$ where the function is $f(x) = 10x + 4$, and do the integration in two segments

$$f_L^{(\alpha)}(x) = \frac{1}{\Gamma(1-\alpha)} \frac{d}{dx} \left( \int_0^{1/2} + \int_{1/2}^x \right) (x-\xi)^{-\alpha} [f(\xi) - f(0)] d\xi$$

$$= \frac{1}{\Gamma(1-\alpha)} \frac{d}{dx} \left( \int_0^{1/2} 10\xi(x-\xi)^{-\alpha} d\xi + \int_{1/2}^x (x-\xi)^{-\alpha} \{17 - 16\xi - 4\} d\xi \right)$$

$$= \frac{1}{\Gamma(1-\alpha)} \frac{d}{dx} \left( \int_0^{1/2} 10(x - (x-\xi))(x-\xi)^{-\alpha} d\xi + \int_{1/2}^x (x-\xi)^{-\alpha} \{13 - 16x + 16(x-\xi)\} d\xi \right)$$

$$= \frac{1}{\Gamma(1-\alpha)} \frac{d}{dx} \left[ 10 \left( \frac{(x-\xi)^{2-\alpha}}{2-\alpha} - x\frac{(x-\xi)^{1-\alpha}}{1-\alpha} \right)_0^{1/2} - \left( 16\frac{(x-\xi)^{2-\alpha}}{2-\alpha} + (13-16x)\frac{(x-\xi)^{1-\alpha}}{1-\alpha} \right)_{1/2}^x \right]$$

$$= \frac{1}{\Gamma(1-\alpha)} \frac{d}{dx} \left[ \begin{array}{l} -10\frac{(x-1/2)^{2-\alpha} - x^{2-\alpha}}{1-\alpha} - 5\frac{(x-1/2)^{1-\alpha}}{1-\alpha} + 10\frac{(x-1/2)^{2-\alpha} - x^{2-\alpha}}{2-\alpha} \\ + 16\frac{(x-1/2)^{1-\alpha} - x^{1-\alpha}}{1-\alpha} - \frac{(x-1/2)^{2-\alpha}}{2-\alpha} + (x-1)\frac{(x-1/2)^{2-\alpha}}{2-\alpha} \\ + \frac{(x-1/2)^{1-\alpha}}{1-\alpha}(13 - 16x) \end{array} \right]$$

$$= \frac{1}{\Gamma(1-\alpha)} \left[ -(x-1/2)^{1-\alpha} \frac{16\alpha + 10}{1-\alpha} + \frac{10 x^{1-\alpha}}{1-\alpha} - (x-1/2)^{1-\alpha}(16x - 8) \right]$$

$$= \frac{\{-26(x-1/2)^{1-\alpha} + 10 x^{1-\alpha}\}}{\Gamma(2-\alpha)}$$

Therefore

$$f_L^{(\alpha)}(x) = \begin{cases} \dfrac{10x^{1-\alpha}}{\Gamma(2-\alpha)} & \text{for} \quad 0 \le x \le 1/2 \\ \dfrac{\{-6(x-1/2)^{1-\alpha} + 10x^{1-\alpha}\}}{\Gamma(2-\alpha)} & \text{for} \quad 1/2 \le x \le 1 \end{cases}$$

Therefore though the function is not differentiable at $x = 1/2$ but the $\alpha$-order derivative at $x = 1/2$ exists and equals to $f_L^{(\alpha)}(1/2) = \dfrac{10(1/2)^{1-\alpha}}{\Gamma(2-\alpha)}$. The graphical presentation of $f_L^{(\alpha)}(x)$ for different values of alpha is shown in the figure-5, from the figure it clear that $f_L^{(\alpha)}(x)$ exists at the non-differentiable point $x = 1/2$.

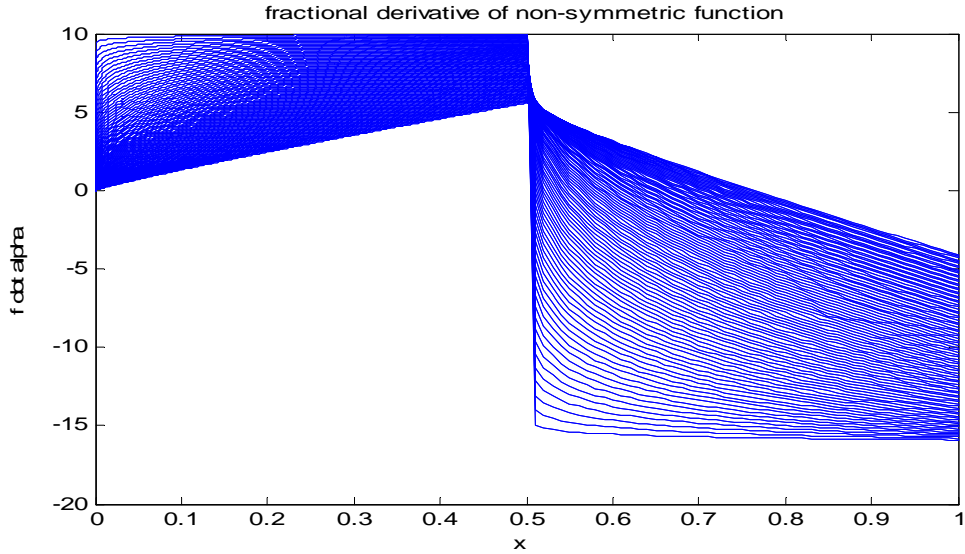

**Fig.5 Graph of the function $f_L^{(\alpha)}(x)$ for different values of alpha.**

The fractional order derivative using right R-L definition on same function we get

$$f_R^{(\alpha)}(x) = -\frac{1}{\Gamma(1-\alpha)} \frac{d}{dx} \int_x^1 (\xi - x)^{-\alpha} [f(1) - f(\xi)] d\xi, \qquad 0 < \alpha < 1.$$

When $0 \le x \le 1/2$

$$f_R^{(\alpha)}(x) = -\frac{1}{\Gamma(1-\alpha)} \frac{d}{dx} \left( \int_x^{1/2} + \int_{1/2}^1 \right) (\xi - x)^{-\alpha} [f(1) - f(\xi)] d\xi$$

$$= \frac{1}{\Gamma(1-\alpha)} \frac{d}{dx} \left( \int_x^{1/2} (\xi - x)^{-\alpha} [10(\xi - x) + (10x + 3)] d\xi + 16 \int_{1/2}^1 [-(\xi - x)^{1-\alpha} + (1-x)(\xi - x)^{-\alpha}] d\xi \right)$$

$$= \frac{1}{\Gamma(1-\alpha)} \frac{d}{dx} \left[ \left[ \frac{10(\xi - x)^{2-\alpha}}{2-\alpha} + (10x+3) \frac{(\xi-x)^{1-\alpha}}{1-\alpha} \right]_x^{1/2} + 16 \left[ (1-x) \frac{(\xi-x)^{1-\alpha}}{1-\alpha} - \frac{(\xi-x)^{2-\alpha}}{2-\alpha} \right]_{1/2}^1 \right]$$

$$= \frac{1}{\Gamma(1-\alpha)} \frac{d}{dx} \left[ \begin{array}{l} \left[ \frac{10(1/2-x)^{2-\alpha}}{2-\alpha} + (10x+3) \frac{(1/2-x)^{1-\alpha}}{1-\alpha} \right] + \\ 16 \left[ (1-x) \frac{(1-x)^{1-\alpha} - (1/2-x)^{1-\alpha}}{1-\alpha} - \frac{(1-x)^{2-\alpha} - (1/2-x)^{2-\alpha}}{2-\alpha} \right] \end{array} \right]$$

$$= \frac{1}{\Gamma(1-\alpha)} \left[ \begin{array}{l} -10(1/2-x)^{1-\alpha} + \frac{10}{1-\alpha}(1/2-x)^{1-\alpha} + 10(1/2-x)^{1-\alpha} \\ -8(1/2-x)^{-\alpha} + 8(1/2-x)^{-\alpha} + \frac{16}{1-\alpha}(1/2-x)^{1-\alpha} - \frac{16}{1-\alpha}(1-x)^{1-\alpha} \end{array} \right]$$

$$= \frac{26(1-x)^{2-\alpha} - 16(1-x)^{1-\alpha}}{\Gamma(2-\alpha)}$$

When $1/2 \leq x \leq 1$

$$f_R^{(\alpha)}(x) = \frac{1}{\Gamma(1-\alpha)} \frac{d}{dx} \int_x^1 (\xi-x)^{-\alpha} (17 - 16\xi - 1) d\xi = \frac{16}{\Gamma(1-\alpha)} \frac{d}{dx} \int_x^1 (\xi-x)^{-\alpha} [(1-x) - (\xi-x)] d\xi$$

$$= \frac{16}{\Gamma(1-\alpha)} \frac{d}{dx} \left[ -\frac{(\xi-x)^{2-\alpha}}{2-\alpha} + (1-x) \frac{(\xi-x)^{1-\alpha}}{1-\alpha} \right]_x^1$$

$$= \frac{16}{\Gamma(1-\alpha)} \frac{d}{dx} \left[ -\frac{(1-x)^{2-\alpha}}{2-\alpha} + (1-x) \frac{(1-x)^{1-\alpha}}{1-\alpha} \right] = -16 \frac{(1-x)^{1-\alpha}}{\Gamma(2-\alpha)}$$

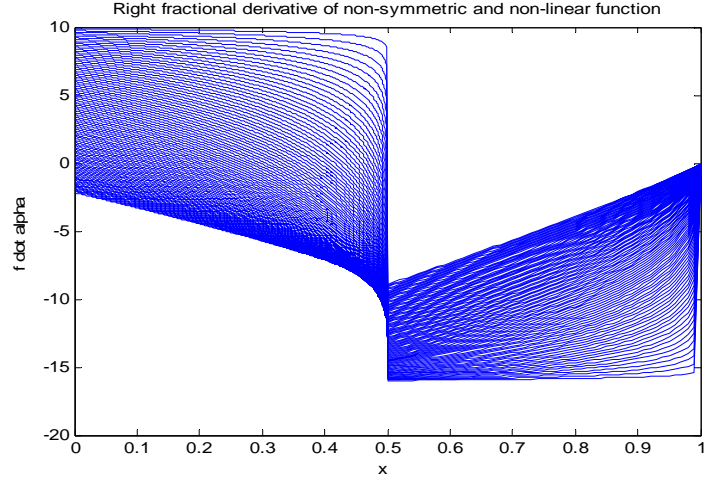

**Fig-6 Graph of the function $f_R^{(\alpha)}(x)$ for different values of alpha.**

Therefore

$$f_R^{(\alpha)}(x) = \begin{cases} \dfrac{26(1-x)^{2-\alpha} - 16(1-x)^{1-\alpha}}{\Gamma(2-\alpha)}, & \text{for } 0 \leq x \leq 1/2 \\ -16\dfrac{(1-x)^{1-\alpha}}{\Gamma(2-\alpha)}, & \text{for } 1/2 \leq x \leq 1 \end{cases}$$

Thus; though the considered function is not differentiable at $x = 1/2$ but its right modified fractional derivative exists and its value is $f_R^{(\alpha)}(1/2) = -16\dfrac{(1/2)^{1-\alpha}}{\Gamma(2-\alpha)}$ which differ from the value $f_L^{(\alpha)}(1/2) = 10\dfrac{(1/2)^{1-\alpha}}{\Gamma(2-\alpha)}$ of the derivative at $x = 1/2$ obtained by left modified R-L derivative by Jumarie modification.

Here the difference indicates there is a phase transition from the left hand to the right hand side about the point $x = 1/2$ and the level or degree of phase transition is $26\dfrac{(1/2)^{1-\alpha}}{\Gamma(2-\alpha)}$.

**Example 3**: Let

$$f(x) = \begin{cases} 30x + 4, & 0 \leq x \leq 0.5 \\ 34 - 30x, & 0.5 \leq x \leq 1 \end{cases}$$

which arises in mapping of a lead in ECG graph. This function is continuous for all values of $x$ in $[0,1]$ but not differentiable at $x = 1/2$.

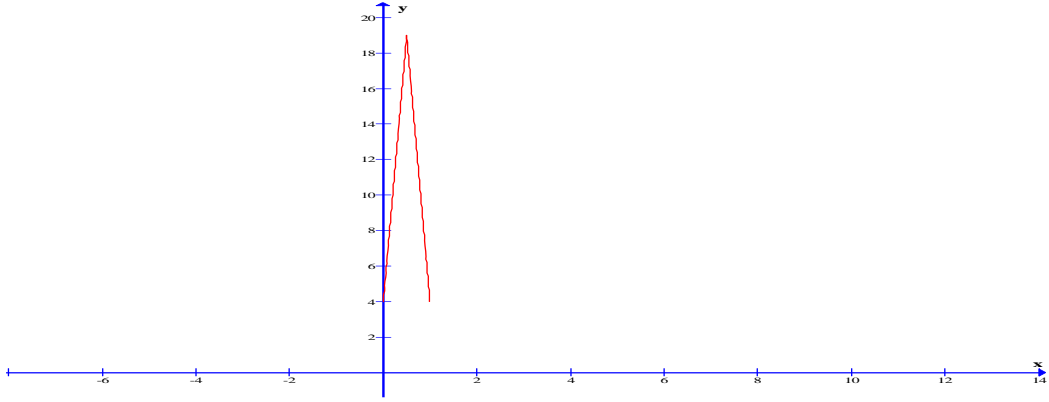

**Fig- 7 Graph of the function** $f(x) = \begin{cases} 30x+4, 0 \leq x \leq 0.5 \\ 34-30x, 0.5 \leq x \leq 1 \end{cases}$

(a) The fractional order derivative using Jumarie modified definition is

$$f_L^{(\alpha)}(x) = \frac{1}{\Gamma(1-\alpha)} \frac{d}{dx} \int_0^x (x-\xi)^{-\alpha} [f(\xi) - f(0)] d\xi, \qquad 0 < \alpha < 1$$

$$= \frac{1}{\Gamma(1-\alpha)} \frac{d}{dx} \int_0^x (x-\xi)^{-\alpha} [f(\xi) - 4] d\xi$$

When $0 \leq x \leq 1/2$

$$f_L^{(\alpha)}(x) = \frac{30}{\Gamma(1-\alpha)} \frac{d}{dx} \int_0^x (x-\xi)^{-\alpha} \xi \, d\xi$$

$$= \frac{30}{\Gamma(1-\alpha)} \frac{d}{dx} \int_0^x (x-\xi)^{-\alpha} \{x - (x-\xi)\} d\xi$$

$$= \frac{30}{\Gamma(1-\alpha)} \frac{d}{dx} \left[ -x \frac{(x-\xi)^{1-\alpha}}{1-\alpha} + \frac{(x-\xi)^{2-\alpha}}{2-\alpha} \right]_0^x = \frac{30}{\Gamma(1-\alpha)} \frac{d}{dx} \left[ \frac{x^{2-\alpha}}{1-\alpha} - \frac{x^{2-\alpha}}{2-\alpha} \right] = \frac{30 x^{1-\alpha}}{\Gamma(2-\alpha)}$$

When $1/2 \leq x \leq 1$.

$$f_L^{(\alpha)}(x) = \frac{30}{\Gamma(1-\alpha)} \frac{d}{dx}\left[\int_0^{1/2} (x-\xi)^{-\alpha}(\xi)d\xi + \int_{1/2}^x (x-\xi)^{-\alpha}(1-\xi)d\xi\right]$$

$$= \frac{30}{\Gamma(1-\alpha)} \frac{d}{dx}\left[\int_0^{1/2} (x-\xi)^{-\alpha}\{x-(x-\xi)\}d\xi + \int_{1/2}^x (x-\xi)^{-\alpha}\{(x-\xi)+(1-x)\}d\xi\right]$$

$$= \frac{30}{\Gamma(1-\alpha)} \frac{d}{dx}\left\{\left[-x\frac{(x-\xi)^{1-\alpha}}{1-\alpha} + \frac{(x-\xi)^{2-\alpha}}{2-\alpha}\right]_0^{1/2} + \left[-\frac{(x-\xi)^{2-\alpha}}{2-\alpha} - (1-x)\frac{(x-\xi)^{1-\alpha}}{1-\alpha}\right]_{1/2}^x\right\}$$

$$= \frac{30}{\Gamma(1-\alpha)} \frac{d}{dx}\left\{\left[-\frac{x(x-\xi)^{1-\alpha} - x^{2-\alpha}}{1-\alpha} + \frac{(x-1/2)^{2-\alpha} - x^{2-\alpha}}{2-\alpha}\right]_0^{1/2} + \left[\frac{(x-1/2)^{2-\alpha}}{2-\alpha} + (1-x)\frac{(x-1/2)^{1-\alpha}}{1-\alpha}\right]_{1/2}^x\right\}$$

$$= \frac{30\{x^{1-\alpha} - 2(x-1/2)^{1-\alpha}\}}{\Gamma(2-\alpha)}$$

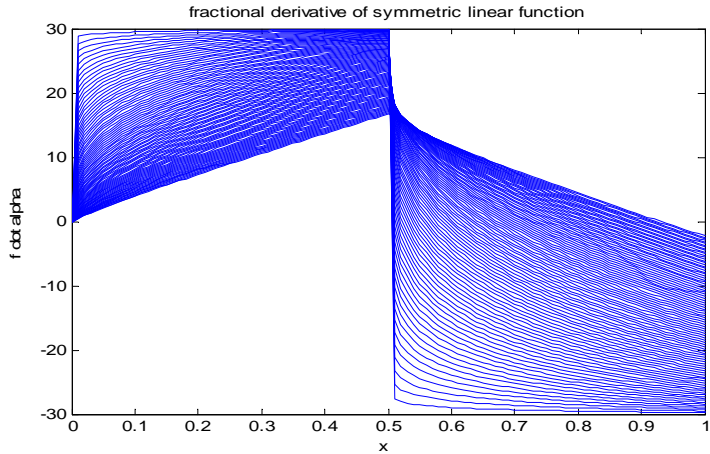

**Fig-8 Graph of the function $f_L^{(\alpha)}(x)$ for different values of alpha.**

From the figure 7 and 8 it is clear that though this function is not differentiable at $x = 1/2$ but its fractional derivative of order $\alpha$ with $0 < \alpha < 1$ exists at $x = 1/2$.

Therefore

$$f_L^{(\alpha)}(x) = \begin{cases} \dfrac{30 x^{1-\alpha}}{\Gamma(2-\alpha)} & \text{for } 0 \leq x \leq 1/2 \\ \dfrac{30\{x^{1-\alpha} - 2(x-1/2)^{1-\alpha}\}}{\Gamma(2-\alpha)} & \text{for } 1/2 \leq x \leq 1 \end{cases}$$

Here

$$f_L^{(\alpha)}(x=0.5) = \frac{30}{\Gamma(2-\alpha)}(1/2)^{1-\alpha}.$$

In previous all the problems we consider the functions which are linear in both side of the non-differentiable point. In the next example we consider a function which is linear in one side and non-linear in other side of the non- differentiable point.

**Example 4**:  Let

$$f(x) = \begin{cases} 4x^2 + 2x + 2, & 0 \leq x \leq 0.5 \\ 5 - 2x, & 0.5 \leq x \leq 1 \end{cases}$$

This function is continuous for all values of $x$ in $[0,1]$ but not differentiable at $x = 1/2$ which is clear from figure-9.

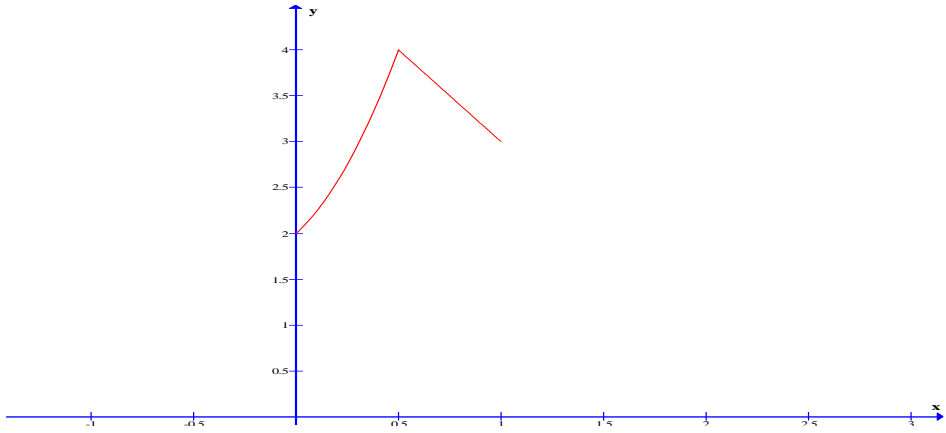

**Fig-9 graph of the function defined above.**

(a) The fractional order derivative using Jumarie modified definition is

$$f_L^{(\alpha)}(x) = \frac{1}{\Gamma(1-\alpha)} \frac{d}{dx} \int_0^x (x-\xi)^{-\alpha}[f(\xi) - f(0)]d\xi, \qquad 0 < \alpha < 1$$

For $0 \leq x \leq 1/2$ $\qquad f(x) = 4x^2 + 2x + 2$

$$f_L^{(\alpha)}(x) = \frac{1}{\Gamma(1-\alpha)} \frac{d}{dx} \int_0^x (x-\xi)^{-\alpha} [4\xi^2 + 2\xi] d\xi$$

$$= \frac{1}{\Gamma(1-\alpha)} \frac{d}{dx} \int_0^x (x-\xi)^{-\alpha} [4\xi^2 + 2\xi + 4x^2 - 8x\xi - 8x^2 + 8x\xi - 2x + 2x + 4x^2] d\xi$$

$$= \frac{1}{\Gamma(1-\alpha)} \frac{d}{dx} \int_0^x (x-\xi)^{-\alpha} [4(x-\xi)^2 - (8x+2)(x-\xi) + 2x + 4x^2] d\xi$$

$$= \frac{1}{\Gamma(1-\alpha)} \frac{d}{dx} \int_0^x \left[ 4(x-\xi)^{2-\alpha} - (8x+2)(x-\xi)^{1-\alpha} + 2x(1+2x^2)(x-\xi)^{-\alpha} \right] d\xi$$

$$= \frac{1}{\Gamma(1-\alpha)} \frac{d}{dx} \left[ -\frac{4(x-\xi)^{3-\alpha}}{3-\alpha} + (8x+2)\frac{(x-\xi)^{2-\alpha}}{2-\alpha} - (2x+4x^2)\frac{(x-\xi)^{1-\alpha}}{1-\alpha} \right]_0^x$$

$$= \frac{1}{\Gamma(1-\alpha)} \frac{d}{dx} \left[ \frac{4x^{3-\alpha}}{3-\alpha} - (8x+2)\frac{x^{2-\alpha}}{2-\alpha} + (2x+4x^2)\frac{x^{1-\alpha}}{1-\alpha} \right]$$

$$= \frac{1}{\Gamma(1-\alpha)} \left[ 2x^{1-\alpha} + 8x^{2-\alpha} \left( \frac{1}{1-\alpha} - \frac{1}{2-\alpha} \right) \right] = \frac{1}{\Gamma(1-\alpha)} \left[ \frac{2x^{1-\alpha}}{1-\alpha} + 8x^{2-\alpha} \left( \frac{1}{1-\alpha} - \frac{1}{2-\alpha} \right) \right]$$

$$= \frac{1}{\Gamma(2-\alpha)} \left[ 2x^{1-\alpha} + \frac{8x^{2-\alpha}}{2-\alpha} \right]$$

$$f_L^{(\alpha)}(1/2) = \frac{(4-\alpha)}{\Gamma(3-\alpha)} \left( \frac{1}{2} \right)^{-\alpha}$$

For $1/2 \leq x \leq 1$ $\quad f(x) = 5 - 2x$

For calculations for the fractional derivative for $1/2 \leq x \leq 1$ we need to calculate, the fractional derivative from start point $f(0)$, and thus take the function in the region $0 \leq x \leq 1/2$ which is $f(x) = 4x^2 + 2x + 2$ and start point value is $f(0) = 2$; and do the integration in two segments first in $[0, 0.5]$ and then $[0.5, 1]$, as demonstrated in the following steps.

$$f_L^{(\alpha)}(x) = \frac{1}{\Gamma(1-\alpha)} \frac{d}{dx} \int_0^x (x-\xi)^{-\alpha}[f(\xi) - f(0)]d\xi, 0 < \alpha < 1 \qquad f(0) = 2$$

$$= \frac{1}{\Gamma(1-\alpha)} \frac{d}{dx} \left[ \int_0^{1/2} (x-\xi)^{-\alpha}[4(x-\xi)^2 - (8x+2)(x-\xi) + 2x + 4x^2]d\xi + \int_{1/2}^x (x-\xi)^{-\alpha}[5 - 2\xi - 2]d\xi \right]$$

$$= \frac{1}{\Gamma(1-\alpha)} \frac{d}{dx} \left[ \begin{array}{l} \int_0^{1/2} (x-\xi)^{-\alpha}[4(x-\xi)^2 - (8x+2)(x-\xi) + 2x + 4x^2]d\xi + \\ \int_{1/2}^x (x-\xi)^{-\alpha}[(3-2x) + 2(x-\xi)]d\xi \end{array} \right]$$

$$= \frac{1}{\Gamma(1-\alpha)} \frac{d}{dx} \left[ \begin{array}{l} \left[ -\frac{4(x-\xi)^{3-\alpha}}{3-\alpha} + (8x+2)\frac{(x-\xi)^{2-\alpha}}{2-\alpha} - (2x+4x^2)\frac{(x-\xi)^{1-\alpha}}{1-\alpha} \right]_0^{1/2} - \\ \left[ (3-2x)\frac{(x-\xi)^{1-\alpha}}{1-\alpha} + 2\frac{(x-\xi)^{2-\alpha}}{2-\alpha} \right]_{1/2}^x \end{array} \right]$$

$$= \frac{1}{\Gamma(1-\alpha)} \frac{d}{dx} \left[ \begin{array}{l} 4\frac{x^{3-\alpha} - (x-1/2)^{3-\alpha}}{3-\alpha} - (8x+2)\frac{x^{2-\alpha} - (x-1/2)^{2-\alpha}}{2-\alpha} + \\ (2x+4x^2)\frac{x^{1-\alpha} - (x-1/2)^{1-\alpha}}{1-\alpha} + 2\frac{(x-1/2)^{1-\alpha}}{1-\alpha} - 2\frac{(x-1/2)^{2-\alpha}}{1-\alpha} + 2\frac{(x-1/2)^{2-\alpha}}{2-\alpha} \end{array} \right]$$

$$= 8\frac{x^{2-\alpha}}{\Gamma(3-\alpha)} + 2\frac{x^{1-\alpha}}{\Gamma(2-\alpha)} - 8\frac{(x-1/2)^{2-\alpha}}{\Gamma(3-\alpha)} - 8\frac{(x-1/2)^{1-\alpha}}{\Gamma(2-\alpha)}$$

Here also

$$f_L^{(\alpha)}(1/2) = \frac{(4-\alpha)}{\Gamma(3-\alpha)} \left(\frac{1}{2}\right)^{-\alpha}$$

$$f_L^{(\alpha)}(x) = \begin{cases} \dfrac{1}{\Gamma(2-\alpha)}\left[2.x^{1-\alpha} + \dfrac{8x^{2-\alpha}}{2-\alpha}\right] & \text{for } 0 \leq x \leq 0.5 \\[2ex] 8\dfrac{x^{2-\alpha}}{\Gamma(3-\alpha)} + 2\dfrac{x^{1-\alpha}}{\Gamma(2-\alpha)} - 8\dfrac{(x-1/2)^{2-\alpha}}{\Gamma(3-\alpha)} - 8\dfrac{(x-1/2)^{1-\alpha}}{\Gamma(2-\alpha)} & \text{for } 0.5 \leq x \leq 1 \end{cases}$$

Thus though this function was non-differentiable at $x = 1/2$ but the fractional derivative exists at $x = 1/2$. The graph of the fractional derivative for different values $\alpha$, $0 < \alpha < 1$ is shown in figure-10.

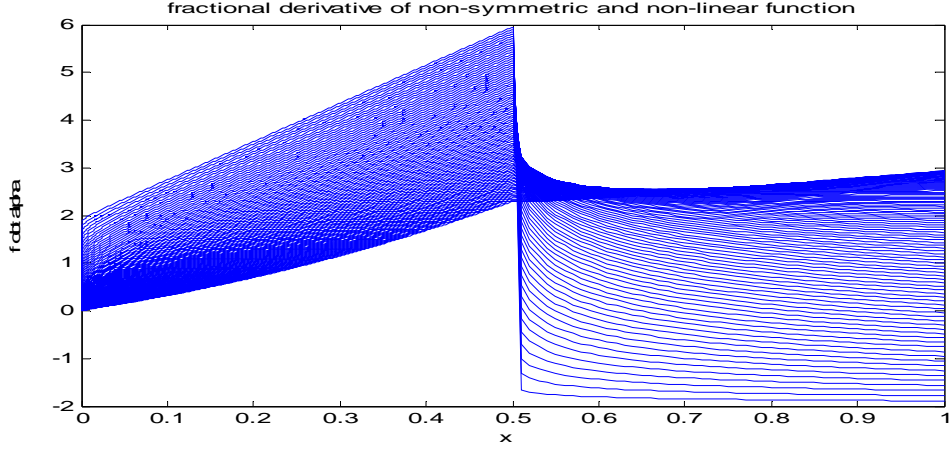

**Fig-10 Graph of the function $f_L^{(\alpha)}(x)$ for different values of alpha.**

(b) The fractional order derivative using right R-L definition on same function we get

$$f_R^{(\alpha)}(x) = -\frac{1}{\Gamma(1-\alpha)} \frac{d}{dx} \int_x^1 (\xi - x)^{-\alpha} [f(1) - f(\xi)] d\xi, \qquad 0 < \alpha < 1$$

When $0 \le x \le 1/2$ $\qquad f(x) = 4x^2 + 2x + 2$

This calculations requires $f(1) = 3$ the end point of the function and the function in the interval $[0.5,1]$ which is $f(x) = 5 - 2x$. We need to carry this integration in two segments as demonstrated below.

$$f_R^{(\alpha)}(x) = -\frac{1}{\Gamma(1-\alpha)} \frac{d}{dx} \left( \int_x^{1/2} + \int_{1/2}^1 \right)(\xi - x)^{-\alpha}[f(1) - f(\xi)]d\xi \qquad f(1) = 3$$

$$= \frac{1}{\Gamma(1-\alpha)} \frac{d}{dx} \begin{pmatrix} \int_x^{1/2} (\xi - x)^{-\alpha}[4(\xi - x)^2 + (8x+2)(\xi - x) + (4x^2 + 2x - 1)]d\xi + \\ 2\int_{1/2}^1 (\xi - x)^{-\alpha}[(1-x) - (\xi - x)]d\xi \end{pmatrix}$$

$$= \frac{1}{\Gamma(1-\alpha)} \frac{d}{dx} \left[ \begin{array}{l} \left[ 4\frac{(\xi-x)^{3-\alpha}}{3-\alpha} + (8x+2)\frac{(\xi-x)^{2-\alpha}}{2-\alpha} + (4x^2+2x-1)\frac{(\xi-x)^{1-\alpha}}{1-\alpha} \right]_x^{1/2} + \\ 2\left[ (1-x)\frac{(\xi-x)^{1-\alpha}}{1-\alpha} - \frac{(\xi-x)^{2-\alpha}}{2-\alpha} \right]_{1/2}^1 \end{array} \right]$$

$$= \frac{1}{\Gamma(1-\alpha)} \frac{d}{dx} \left[ \begin{array}{l} \left[ 4\frac{(1/2-x)^{3-\alpha}}{3-\alpha} + (8x+2)\frac{(1/2-x)^{2-\alpha}}{2-\alpha} + (4x^2+2x-1)\frac{(1/2-x)^{1-\alpha}}{1-\alpha} \right] + \\ 2\left[ (1-x)\frac{(1-x)^{1-\alpha}-(1/2-x)^{1-\alpha}}{1-\alpha} - \frac{(1-x)^{2-\alpha}-(1/2-x)^{2-\alpha}}{2-\alpha} \right] \end{array} \right]$$

$$= -\frac{2(1-x)^{1-\alpha}}{\Gamma(2-\alpha)} + \frac{8(1/2-x)^{1-\alpha}}{\Gamma(2-\alpha)} - \frac{8(1-x)^{2-\alpha}}{\Gamma(3-\alpha)}$$

When $1/2 \leq x \leq 1$ $\quad f(x) = 5 - 2x \quad\quad f(1) = 3$

$$f_R^{(\alpha)}(x) = \frac{1}{\Gamma(1-\alpha)} \frac{d}{dx} \int_x^1 (\xi - x)^{-\alpha} (5 - 2\xi - 3)] d\xi$$

$$= \frac{2}{\Gamma(1-\alpha)} \frac{d}{dx} \int_x^1 (\xi - x)^{-\alpha} [(1-x) - (\xi - x)] d\xi$$

$$= \frac{2}{\Gamma(1-\alpha)} \frac{d}{dx} \left[ (1-x)\frac{(\xi-x)^{1-\alpha}}{1-\alpha} - \frac{(\xi-x)^{2-\alpha}}{2-\alpha} \right]_x^1$$

$$= \frac{16}{\Gamma(1-\alpha)} \frac{d}{dx} \left[ \frac{(1-x)^{2-\alpha}}{1-\alpha} - \frac{(1-x)^{2-\alpha}}{2-\alpha} \right] = -\frac{2(1-x)^{1-\alpha}}{\Gamma(2-\alpha)}$$

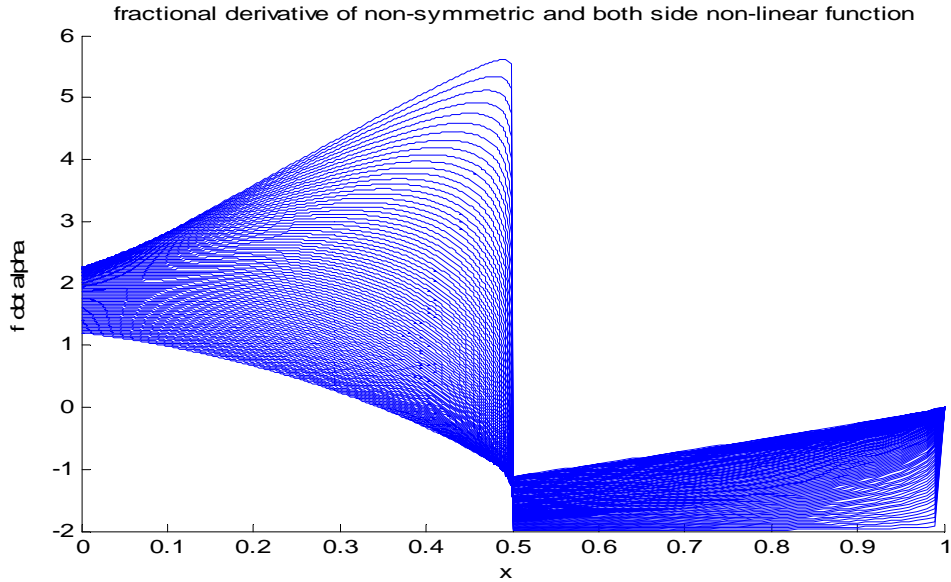

**Fig-11: Graph of the function $f_R^{(\alpha)}(x)$ for different values of alpha.**

Therefore

$$f_R^{(\alpha)}(x) = \begin{cases} -\dfrac{2(1-x)^{1-\alpha}}{\Gamma(2-\alpha)} + \dfrac{8(1/2-x)^{1-\alpha}}{\Gamma(2-\alpha)} - \dfrac{8(1-x)^{2-\alpha}}{\Gamma(3-\alpha)}, & \text{for } 0 \leq x \leq 1/2 \\ -\dfrac{2(1-x)^{1-\alpha}}{\Gamma(2-\alpha)}, & \text{for } 1/2 \leq x \leq 1 \end{cases}$$

The fractional right modified derivative of the function is shown in the figure-11. The value of the derivative at $x = 1/2$ is

$$f_R^{(\alpha)}(1/2) = \frac{-2}{\Gamma(3-\alpha)} \left(\frac{1}{2}\right)^{1-\alpha}$$

This value of right derivative differs from the value of the Jumarie left derivative to this function which indicates there is a phase transition at $x = 1/2$.

The previous function was linear in one side of the differentiable function and non-linear in other side, of the transition point. Now we consider a function which is non-linear in both sides with respect to the non-differentiable point, but continuous at that transition point.

**Example 5**: Now consider the function

$$f(x) = \begin{cases} 4x^2 + 3, & 0 \leq x \leq 0.5 \\ 5 - 4x^2, & 0.5 \leq x \leq 1 \end{cases}.$$

This function is not symmetric with respect to the non-differentiable point $x = 1/2$, refer figure-12. To characterize this function the fractional derivative of this function is calculated using both side modified definition of derivative.

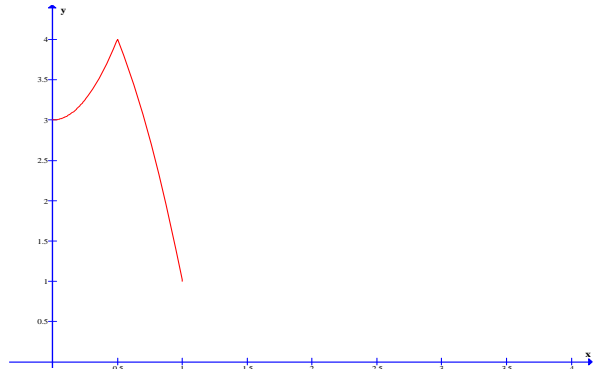

**Fig-12 graph of the function defined above.**

(a) The fractional order derivative using Jumarie modified definition is

$$f_L^{(\alpha)}(x) = \frac{1}{\Gamma(1-\alpha)} \frac{d}{dx} \int_0^x (x-\xi)^{-\alpha} [f(\xi) - f(0)] d\xi, \qquad 0 < \alpha < 1$$

For $0 \leq x \leq 0.5$ $\quad f(x) = 4x^2 + 3$ $\quad f(0) = 3$

$$f_L^{(\alpha)}(x) = \frac{1}{\Gamma(1-\alpha)} \frac{d}{dx} \int_0^x (x-\xi)^{-\alpha} [f(\xi) - f(0)] d\xi$$

$$= \frac{1}{\Gamma(1-\alpha)} \frac{d}{dx} \int_0^x (x-\xi)^{-\alpha} [4(x-\xi)^2 - 8x(x-\xi) + 4x^2] d\xi$$

$$= \frac{1}{\Gamma(1-\alpha)} \frac{d}{dx} \left[ -\frac{4(x-\xi)^{3-\alpha}}{3-\alpha} + \frac{8x(x-\xi)^{2-\alpha}}{2-\alpha} - 4x^2 \frac{(x-\xi)^{1-\alpha}}{1-\alpha} \right]_0^x$$

$$= \frac{1}{\Gamma(1-\alpha)} \frac{d}{dx} \left[ \frac{4x^{3-\alpha}}{3-\alpha} + \frac{8x^{3-\alpha}}{2-\alpha} - 4\frac{x^{3-\alpha}}{1-\alpha} \right] = \frac{8x^{2-\alpha}}{\Gamma(3-\alpha)}$$

For $0.5 \leq x \leq 1$ $\quad f(x) = 5 - 4x^2$

In this calculations we require the value at the start point which is $f(0) = 3$ also requiring the function in the interval $[0, 0.5]$ which is $f(x) = 4x^2 + 3$, and integration is done in two segments $[0, 0.5]$ and $[0.5, 1]$

$$f_L^{(\alpha)}(x) = \frac{1}{\Gamma(1-\alpha)} \frac{d}{dx} \int_0^x (x-\xi)^{-\alpha} [f(\xi) - f(0)] d\xi$$

$$= \frac{1}{\Gamma(1-\alpha)} \frac{d}{dx} \left[ \int_0^{1/2} (x-\xi)^{-\alpha} \left[ 4\xi^2 + 3 - 3 \right] d\xi + \int_{1/2}^x \left[ 5 - 4\xi^2 - 3 \right] d\xi \right]$$

$$= \frac{1}{\Gamma(1-\alpha)} \frac{d}{dx} \left[ \begin{array}{l} \int_0^{1/2} (x-\xi)^{-\alpha} [4(x-\xi)^2 - 8x(x-\xi) + 4x^2] d\xi + \\ \int_{1/2}^x (x-\xi)^{-\alpha} [(2-4x^2) - 4(x-\xi)^2 + 8x(x-\xi)] d\xi \end{array} \right]$$

$$= \frac{1}{\Gamma(1-\alpha)} \frac{d}{dx} \left[ \begin{array}{l} \left[ -\frac{4(x-\xi)^{3-\alpha}}{3-\alpha} + \frac{8x(x-\xi)^{2-\alpha}}{2-\alpha} - 4x^2 \frac{(x-\xi)^{1-\alpha}}{1-\alpha} \right]_0^{1/2} + \\ \left[ \frac{4(x-\xi)^{3-\alpha}}{3-\alpha} - \frac{8x(x-\xi)^{2-\alpha}}{2-\alpha} - (2-4x^2) \frac{(x-\xi)^{1-\alpha}}{1-\alpha} \right]_{1/2}^x \end{array} \right]$$

$$= \frac{1}{\Gamma(1-\alpha)} \frac{d}{dx} \left[ \begin{array}{l} -4\frac{(x-1/2)^{3-\alpha} - x^{3-\alpha}}{3-\alpha} + \frac{8x\{(x-1/2)^{2-\alpha} - x^{2-\alpha}\}}{2-\alpha} - 4x^2 \frac{(x-1/2)^{1-\alpha} - x^{1-\alpha}}{1-\alpha} - \\ \frac{4(x-1/2)^{3-\alpha}}{3-\alpha} + \frac{8x(x-1/2)^{2-\alpha}}{2-\alpha} + (2-4x^2) \frac{(x-1/2)^{1-\alpha}}{1-\alpha} \end{array} \right]$$

$$= \frac{8x^{2-\alpha}}{\Gamma(3-\alpha)} - \frac{8(x-1/2)^{1-\alpha}}{\Gamma(2-\alpha)} - \frac{16(x-1/2)^{2-\alpha}(\alpha^2 - 3\alpha + 3)}{\Gamma(3-\alpha)}$$

Therefore

$$f_L^{(\alpha)}(x) = \begin{cases} \dfrac{8x^{2-\alpha}}{\Gamma(3-\alpha)}, & 0 \leq x \leq 0.5 \\ \dfrac{8x^{2-\alpha}}{\Gamma(3-\alpha)} - \dfrac{8(x-1/2)^{1-\alpha}}{\Gamma(2-\alpha)} - \dfrac{16(x-1/2)^{2-\alpha}}{\Gamma(3-\alpha)}, & 0.5 \leq x \leq 1 \end{cases}$$

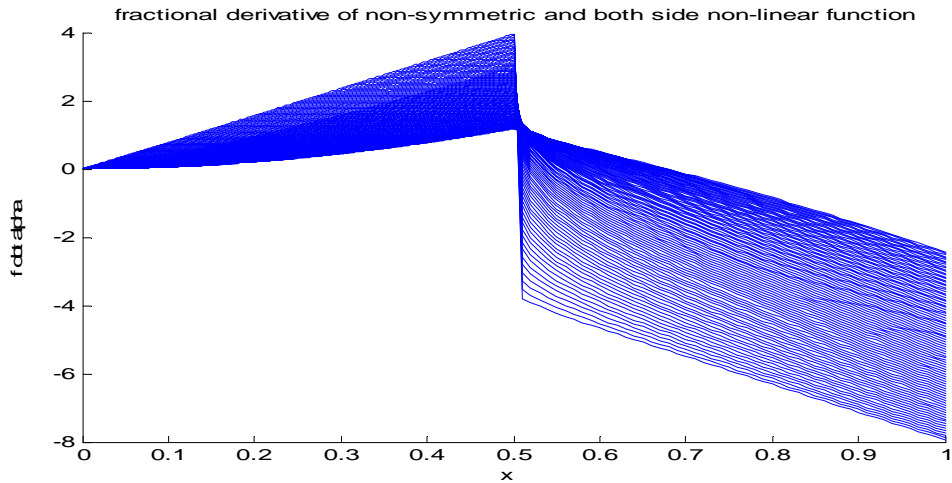

**Fig-13 Graph of the function $f_L^{(\alpha)}(x)$ for different values of alpha.**

(b) The fractional order derivative using right R-L definition on same function we get

$$f(x) = \begin{cases} 4x^2 + 3, 0 \leq x \leq 0.5 \\ 5 - 4x^2, 0.5 \leq x \leq 1 \end{cases}$$

$$f_R^{(\alpha)}(x) = -\frac{1}{\Gamma(1-\alpha)} \frac{d}{dx} \int_x^1 (\xi - x)^{-\alpha} [f(1) - f(\xi)] d\xi, 0 < \alpha < 1.$$

When $0 \leq x \leq 1/2$ $\quad f(x) = 4x^2 + 3 \quad f(1) = 1$

We need here the integration process in two intervals $[0, 0.5]$ and $[0.5, 1]$ as demonstrated below

$$f_R^{(\alpha)}(x) = -\frac{1}{\Gamma(1-\alpha)} \frac{d}{dx} \left( \int_x^{1/2} + \int_{1/2}^1 \right) (\xi-x)^{-\alpha} [f(1) - f(\xi)] d\xi$$

$$= \frac{-1}{\Gamma(1-\alpha)} \frac{d}{dx} \left( \int_x^{1/2} (\xi-x)^{-\alpha}[1-4\xi^2-3]d\xi + \int_{1/2}^1 (\xi-x)^{-\alpha}[1-5+4\xi^2]d\xi \right)$$

$$= \frac{1}{\Gamma(1-\alpha)} \frac{d}{dx} \left( \begin{array}{l} \int_x^{1/2} (\xi-x)^{-\alpha}[4(\xi-x)^2 + 8x(\xi-x) + (4x^2+2)]d\xi + \\ 4\int_{1/2}^1 (\xi-x)^{-\alpha}[(1-x^2) - (\xi-x)^2 - 2x(\xi-x)]d\xi \end{array} \right)$$

$$= \frac{1}{\Gamma(1-\alpha)} \frac{d}{dx} \left[ \begin{array}{l} \left[ 4\frac{(\xi-x)^{3-\alpha}}{3-\alpha} + 8x\frac{(\xi-x)^{2-\alpha}}{2-\alpha} + (4x^2+2)\frac{(\xi-x)^{1-\alpha}}{1-\alpha} \right]_x^{1/2} + \\ 4\left[ (1-x^2)\frac{(\xi-x)^{1-\alpha}}{1-\alpha} - 2x\frac{(\xi-x)^{2-\alpha}}{2-\alpha} - \frac{(\xi-x)^{3-\alpha}}{3-\alpha} \right]_{1/2}^1 \end{array} \right]$$

$$= \frac{1}{\Gamma(1-\alpha)} \frac{d}{dx} \left[ \begin{array}{l} \left[ 4\frac{(1/2-x)^{3-\alpha}}{3-\alpha} + 8x\frac{(1/2-x)^{2-\alpha}}{2-\alpha} + (4x^2+2)\frac{(1/2-x)^{1-\alpha}}{1-\alpha} \right] + \\ 4\left[ (1-x^2)\frac{(1-x)^{1-\alpha} - (1/2-x)^{1-\alpha}}{1-\alpha} - 8x\frac{(1-x)^{2-\alpha} - (1/2-x)^{2-\alpha}}{2-\alpha} - \frac{(1-x)^{3-\alpha} - (1/2-x)^{3-\alpha}}{3-\alpha} \right] \end{array} \right]$$

$$= 8\frac{(1-x)^{2-\alpha}}{\Gamma(3-\alpha)} - 8\frac{(1-x)^{1-\alpha}}{\Gamma(2-\alpha)} + 8\frac{(1/2-x)^{1-\alpha}}{\Gamma(2-\alpha)} - 16\frac{(1-x)^{2-\alpha}}{\Gamma(3-\alpha)}$$

When $1/2 \leq x \leq 1$     $f(x) = 5 - 4x^2$     $f(1) = 1$

$$f_R^{(\alpha)}(x) = \frac{1}{\Gamma(1-\alpha)} \frac{d}{dx} \int_x^1 (\xi-x)^{-\alpha}(5-4\xi^2-1)]d\xi$$

$$= \frac{4}{\Gamma(1-\alpha)} \frac{d}{dx} \int_x^1 (\xi-x)^{-\alpha}[(1-x^2) - (\xi-x)^2 - 2x(\xi-x)]d\xi$$

$$= \frac{16}{\Gamma(1-\alpha)} \frac{d}{dx} \left[ (1-x^2)\frac{(\xi-x)^{1-\alpha}}{1-\alpha} - 2x\frac{(\xi-x)^{2-\alpha}}{2-\alpha} - \frac{(\xi-x)^{3-\alpha}}{3-\alpha} \right]_x^1$$

$$= \frac{16}{\Gamma(1-\alpha)} \frac{d}{dx} \left[ (1-x^2)\frac{(1-x)^{1-\alpha}}{1-\alpha} - 2x\frac{(1-x)^{2-\alpha}}{2-\alpha} - \frac{(1-x)^{3-\alpha}}{3-\alpha} \right] = 8\frac{(1-x)^{2-\alpha}}{\Gamma(3-\alpha)} - 8\frac{(1-x)^{1-\alpha}}{\Gamma(2-\alpha)}$$

Therefore

$$f_R^{(\alpha)}(x) = \begin{cases} 8\dfrac{(1-x)^{2-\alpha}}{\Gamma(3-\alpha)} - 8\dfrac{(1-x)^{1-\alpha}}{\Gamma(2-\alpha)} + 8\dfrac{(1/2-x)^{1-\alpha}}{\Gamma(2-\alpha)} - 16\dfrac{(1-x)^{2-\alpha}}{\Gamma(3-\alpha)}, & \text{for } 0 \leq x \leq 1/2 \\ 8\dfrac{(1-x)^{2-\alpha}}{\Gamma(3-\alpha)} - 8\dfrac{(1-x)^{1-\alpha}}{\Gamma(2-\alpha)}, & \text{for } 1/2 \leq x \leq 1 \end{cases}$$

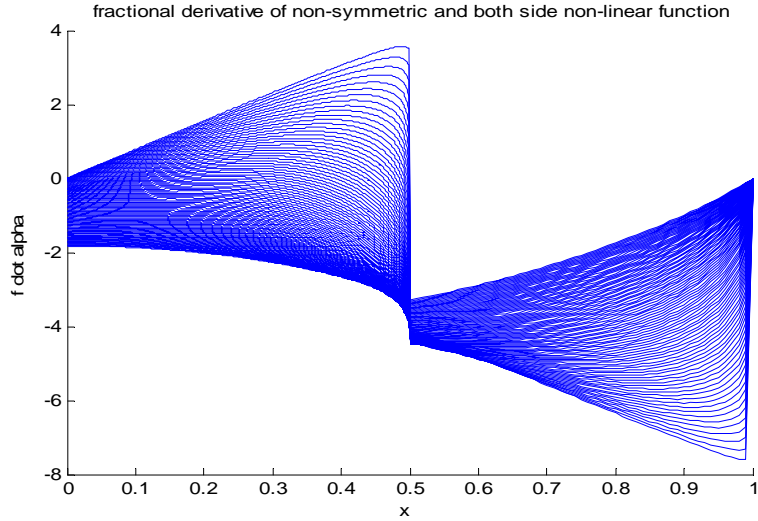

**Fig-14 Graph of the function $f_R^{(\alpha)}(x)$ for different values of alpha.**

### 3.0 Conclusion

From the above examples it is clear that some functions which are unreachable at a non-differentiable points in the defined interval in classical derivative sense but they are differentiable in fractional sense. But the modified fractional derivative in both left and right sense gives same value for differentiable functions but gives different value for non-differentiable case. For non-differentiable the difference in values of the fractional derivative in Jumarie modified and right R-L modified sense indicates there is phase transition about the non-differentiable points. The difference in values indicates the level of the phase transition. These are useful indicators to quantify and compare the non-differentiable but continuous points in a system. This method we are extending to differentiate various ECG graphs by quantification of non-differentiable points; is useful method in differential diagnostic.

### 4.0 Reference

[1] Ross, B., "The development of fractional calculus 1695-1900". *Historia Mathematica,* 1977. 4(1), pp. 75–89.

[2] Diethelm. K. "The analysis of Fractional Differential equations. Springer-Verlag. 2010.

[3] Kilbas A, Srivastava HM, Trujillo JJ. Theory and Applications of Fractional Differential Equations. North-Holland Mathematics Studies, Elsevier Science, Amsterdam, the Netherlands. 2006;204:1-523.

[4] Miller KS, Ross B. An Introduction to the Fractional Calculus and Fractional Differential Equations.John Wiley & Sons, New York, NY, USA; 1993.

[5] Samko SG, Kilbas AA, Marichev OI. Fractional Integrals and Derivatives. Gordon and Breach Science, Yverdon, Switzerland; 1993.

[6] Podlubny I. Fractional Differential Equations, Mathematics in Science and Engineering, Academic Press, San Diego, Calif, USA. 1999;198.

[7] Oldham KB, Spanier J. The Fractional Calculus, Academic Press, New York, NY, USA; 1974.

[8] Das. S. Functional Fractional Calculus 2$^{nd}$ Edition, Springer-Verlag 2011.

[9] Gemant, A., "On fractional differentials". *Phil. Mag. (Ser. 7),* 1942. **25**(4), pp. 540–549.

[10] G. Jumarie, Fractional partial differential equations and modified Riemann-Liouville derivatives. Method for solution. J. Appl. Math. and Computing, 2007 (24), Nos 1-2, pp 31-48.

[11] G. Jumarie, Modified Riemann-Liouville derivative and fractional Taylor series of non-differentiable functions Further results, Computers and Mathematics with Applications, 2006. (51), 1367-1376.

[12] H. Kober, On fractional integrals and derivatives, Quart. J. Math. Oxford, 1940 (11), pp 193-215.

[13] A.V. Letnivov, Theory of differentiation of fractional order, Math. Sb., (3)1868, pp 1-7.

[14] J. Liouville, Sur le calcul des differentielles `a indices quelconques(in french), J. Ecole Polytechnique, (13)1832, 71.

[15] M. Caputo, "Linear models of dissipation whose *q* is almost frequency independent-ii," *Geophysical Journal of the Royal Astronomical Society*, , 1967. vol. 13, no. 5, pp. 529–539.

[16] Abhay Parvate, A. D. Gangal. Calculus of fractals subset of real line: formulation-1; World Scientific, Fractals Vol. 17, 2009.

[17] Abhay Parvate, Seema satin and A.D.Gangal. Calculus on a fractal curve in R$^n$ arXiv:00906 oo76v1 3.6.2009; also in Pramana-J-Phys.

[18] Abhay parvate, A. D. Gangal. Fractal differential equation and fractal time dynamic systems, Pramana-J-Phys, Vol 64, No. 3, 2005 pp 389-409.

[19]. E. Satin, Abhay Parvate, A. D. Gangal. Fokker-Plank Equation on Fractal Curves, Seema, Chaos Solitons & Fractals-52 2013, pp 30-35.

[20] K M Kolwankar and A D. Gangal. Local fractional Fokker plank equation, Phys Rev Lett. 80 1998,

**Acknowledgement**

Acknowledgments are to **Board of Research in Nuclear Science** (BRNS), Department of Atomic Energy Government of India for financial assistance received through BRNS research project no. 37(3)/14/46/2014-BRNS with BSC BRNS, title "Characterization of unreachable (Holderian) functions via Local Fractional Derivative and Deviation Function.